\documentclass[12pt,oneside,frenchb,english,reqno]{amsart}
\usepackage{times}
\usepackage[T1]{fontenc}
\usepackage[latin1]{inputenc}
\usepackage{amssymb}

\makeatletter

\providecommand{\LyX}{L\kern-.1667em\lower.25em\hbox{Y}\kern-.125emX\@}

 \theoremstyle{remark}    
 \newtheorem{claim}{Claim}
 \theoremstyle{remark}
 \newtheorem*{rem*}{Remark}
 \theoremstyle{plain}    
 \newtheorem{thm}{Theorem} 

\def\th@nopoint{
  \thm@headpunct{} %This is what removes the blasted point!
  \itshape % This makes the body in italics, just like a theorem.
}
\theoremstyle{nopoint}
\newtheorem*{ialt}{}

\usepackage{babel}
\makeatother
\begin{document}

\title{A null series with small anti-analytic part.}

\author{Gady Kozma}

\email{gadyk@wisdom.weizmann.ac.il \\
gadykozma@hotmail.com}

\curraddr{The Weizmann Institute of Science, Rehovot, Israel.}

\author{Alexander Olevski\v i}

\email{olevskii@post.tau.ac.il}\address{School of Mathematical Sciences, Tel Aviv University, Ramat Aviv 69978, Israel}

\thanks{Research supported in part by the Israel Science Foundation}

\begin{abstract}
We show that it is possible for an $L^{2}$ function on the circle,
which is a sum of an almost everywhere convergent series of exponentials
with positive frequencies, to not belong to the Hardy space $H^{2}$.
A consequence in the uniqueness theory is obtained.
\end{abstract}
\maketitle

\section{INTRODUCTION}

A (nontrivial) trigonometric series\begin{equation}
\sum c(n)e^{int}\quad t\in \mathbb{T}=\mathbb{R}/2\pi \mathbb{Z}\label{eq:defser}\end{equation}
is called a null series if it converges to zero almost everywhere
(a.e.). The existence of such a series was discovered by D. E. Menshov
in 1916 (see \cite[chap.\ XIV]{B64}). He constructed a singular compactly
supported finite Borel measure on $\mathbb{T}$ with Fourier transform
vanishing at infinity. The Riemannian theory implies that the Fourier
series of this measure converges to zero at every point outside of
the support. This famous example of Menshov was the origin of the
modern uniqueness theory in Fourier Analysis, see \cite{B64,KS94,KL87}.

Clearly a null series can not belong to $L^{2}$. A less trivial observation
is that it can not be {}``analytic'' that is involve positive frequencies
only. This follows from the Abel summation and Privalov {}``angular
limit'' theorems. It turns out however that the {}``non-analytic''
part of a null series may belong to $L^{2}$.

\begin{thm}
\label{thm:Mensh}There is a null series (\ref{eq:defser}) such that\[
\sum _{n<0}|c(n)|^{2}<\infty \quad .\]

\end{thm}
An equivalent formulation of the result :

\begin{thm}
\label{thm:complx}There is a power series\[
F(z)=\sum c(n)z^{n}\]
converging a.e.~on the circle $|z|=1$ to some function $f\in L^{2}(\mathbb{T})$,
and $f$ does not belong to the Hardy space $H^{2}$.
\end{thm}
When we say that a function $f$ on the circle is in $H^{2}$ we mean
that it is a boundary limit of an $H^{2}$ function on the disk, or,
equivalently, that $f \in L^2$ and $\widehat{f}(-n)=0$ for $n=1,2,\dotsc $

To see that theorem \ref{thm:Mensh} implies \ref{thm:complx}, use
Carleson's convergence theorem \cite{C66} to get that the analytic
part of the sum, $\sum _{n\geq 0}c(n)z^{n}$ converges a.e.~on the
circle $|z|=1$ and then use, as above, Abel and Privalov's
theorems to get that the resulting $f$ is not in $H^{2}$.
Reversing these arguments one may derive theorem \ref{thm:Mensh}
from theorem \ref{thm:complx}.

It should be mentioned that usually if representation by harmonics
is unique then it is the Fourier series. Compare for instance the
classical Cantor and du Bois-Reymond \cite[pp.\ 193, 201]{B64} theorems
on pointwise convergence everywhere. Our result shows that this principle
is not universal. Indeed, any $f$ may have \emph{at most one representation}
by an a.e.~convergent series \begin{equation}
\sum _{n\geq 0}c(n)e^{int}\label{eq:sumanal}\end{equation}
however, even if $f\in L^{2}$, the coefficients, in general, can
not be recovered by Fourier's formula.

\section{PROOF}

\subsection{\label{sub:defFf}}

Our main goal is to construct an {}``analytic pseudofunction'',
that is\begin{align}
F(z) & =\sum _{n\geq 0}c(n)z^{n}\label{eq:defF}\\
c(n) & =o(1)\label{eq:cno1}
\end{align}
with the following properties:

\begin{enumerate}
\item \label{enu:FnotinH2}$F\not \in H^{2}$.
\item There is a compact $K\subset \mathbb{T}$ of Lebesgue measure zero
such that $F$ has boundary values on $K^{c}$\[
f(t):=\lim _{z\rightarrow e^{it}}F(z)\quad \forall t\not \in K.\]

\item $f\in L^{\infty }(\mathbb{T})$.
\item \label{enu:limunif}The limit is uniform on any closed arc $J\subset K^{c}$.
\end{enumerate}
Having a function $F$ with all the properties above one can get the
result easily. Indeed, the series (\ref{eq:defF}) on the boundary
represents a distribution\[
\bar{F}:=\sum \widehat{F}(n)e^{int},\]
which is the limit (in distributional sense) of $F_{r}:=F(re^{it})$
as $r \rightarrow 1$.
On the other hand $F(re^{it})\rightarrow f(t)$ uniformly on any closed arc
$J\subset K^{c}$, so the distribution $\bar{F}-f$ is supported on
$K$. Hence the condition (\ref{eq:cno1}) implies uniform convergence
of the Fourier series of $\bar{F}-f$ to zero on any such $J$, see
\cite[p.\ 54]{KS94}. Theorem \ref{thm:complx} (and hence theorem
\ref{thm:Mensh}) will follow.

The function $F$ will be obtained as $1/G$ where $G$ is a singular
inner function, so\begin{equation}
F(z)=\exp \Big (\int _{\mathbb{T}}\frac{e^{it}+z}{e^{it}-z}d\mu (t)\Big ).\label{eq:poisson}\end{equation}
This construction will ensure (\ref{enu:FnotinH2})-(\ref{enu:limunif}),
if $\mu $ is a positive measure supported on $K$, so our task in
sections \ref{sub:lnKnun}-\ref{sub:defmu} will be to construct a
singular $\mu $ such that (\ref{eq:cno1}) will be satisfied.

\subsection{\label{sub:lnKnun}}

Denoting $g(x):=xe^{2/x}+1-x$ we fix a sequence\[
l(1)>l(2)>...\rightarrow 0,\]
such that\[
g(l(n))-g(l(n-1))=o(1).\]
Proceed with the induction as follows. Let $K_{0}=\mathbb{T}$. Suppose
we already have a compact $K_{n-1}\subset \mathbb{T}$ which is a
finite union of segments of equal lengths. Divide each of them to
$q(n)$ equal subsegments $I$ and replace each $I$ by the concentric
segment $I'$,\[
|I'|=\frac{l(n)}{l(n-1)}|I|\]
(here and below by $|E|$ we denote the normalized Lebesgue measure
of a set $E\subset \mathbb{T}$). Set $K_{n}:=\cup I'$, so $|K_{n}|=l(n)$,
and\[
u_{n}:=\frac{1}{l(n)}\mathbf{1}_{K_{n}}.\]

\begin{claim}
If the number $q=q(n)$ is sufficiently large then the function $u_{n}-u_{n-1}$
is {}``almost orthogonal'' to any pre-given finite dimensional subspace
in $L^{2}(\mathbb{T})$. More precisely: for any $\epsilon >0$, $N\in \mathbb{N}$
there is a $Q\in \mathbb{N}$ such that $\forall q(n)>Q$,\[
|\widehat{u_{n}-u_{n-1}}(k)|<\epsilon \quad \forall k,|k|<N\]
(here and below the sign $\widehat{\cdot }$ stands for the Fourier
transform on $\mathbb{T}$. The term {}``sufficiently large'' means
that the minimal allowed value may depend on everything that happened
in previous stages of the induction). 
\end{claim}
To prove the claim it is enough to mention that $u_{n}-u_{n-1}$ is
supported on the union of the segments $I$, the length of each $I$
is arbitrary small as $q$ gets large and the average of $u_{n}-u_{n-1}$
on $I$ equals to zero.\qed

Obviously the same inequality holds for the conjugate function $\widetilde{u_{n}-u_{n-1}}$
so we obtain

\begin{claim}
\label{cla:ununtilde}Given $\epsilon >0$, $N\in \mathbb{N}$ and
sufficiently large $q$ the function\[
h_{n}:=(u_{n}-u_{n-1})+i(\widetilde{u_{n}-u_{n-1}})\]
satisfies \[
|\widehat{h_{n}}(k)|<\epsilon \quad \forall k,0\leq k<N\]

\end{claim}
Now denote: $f_{n}=e^{u_{n}+i\widetilde{u_{n}}}$.

\begin{claim}
\label{cla:fnfn1korth}Given $\epsilon >0$, $N\in \mathbb{N}$ and
$q$ sufficiently large we have :\[
|\widehat{f_{n}-f_{n-1}}(k)|<\epsilon \quad \forall k,0\leq k<N\]

\end{claim}
Indeed, \[
f_{n}-f_{n-1}=f_{n-1}(e^{h_{n}}-1).\]
Clearly the fact that $h_{n}$ are analytic gives that they may be
exponentiated formally (e.g.~by extending to the disk $\mathbb{D}$
and using $\widehat{h}(k)=h^{(k)}(0)$) which gives that $\widehat{e^{h_{n}}-1}(k)$
is a polynomial with no constant term in $\widehat{h_{n}}(1),\dotsc ,\linebreak [0]\widehat{h_{n}}(k)$.
Since $f_{n-1}$ is also analytic, $\widehat{f_{n}-f_{n-1}}(k)$ is
a finite combination of $\widehat{f_{n-1}}(j)$ and $\widehat{e^{h_{n}}-1}(k-j)$
and the claim is a consequence of claim \ref{cla:ununtilde}.

As an immediate corollary we obtain :

\begin{claim}
\label{cla:fnfn1fn1orth}For any $\epsilon >0$ and sufficiently large
$q$\[
|\langle f_{n}-f_{n-1},f_{n-1}\rangle |<\epsilon .\]

\end{claim}
We mean here the usual inner product, $\langle f,g\rangle =\int _{\mathbb{T}}f\bar{g}$.

\subsection{}

Proceeding with the induction above we get

\begin{claim}
If the numbers $q(n)$ grow sufficiently fast then there are numbers
$N_{1}<N_{2}<\dotsm $ such that the functions $\{f_{n}\}$ satisfy,
for any $n$, the conditions:
\end{claim}
\begin{enumerate}
\item \label{enu:defNn}$|\widehat{f_{n-1}}(k)|<\frac{1}{n}$, for all $k$
such that $k>N_{n}$.
\item $|\widehat{f_{n}-f_{n-1}}(k)|<2^{-n}$, for all $k$ such that $0\leq k\leq N_{n}$.\label{enu:fnfn1k}
\item $|\langle f_{n}-f_{n-1},f_{n}\rangle |<\frac{1}{n}$.\label{enu:fnorth}
\end{enumerate}
It is enough on the $n$th step of the induction to choose $N_{n}$
so that (\ref{enu:defNn}) is fulfilled and then to use claims \ref{cla:fnfn1korth}
and \ref{cla:fnfn1fn1orth} to ensure (\ref{enu:fnfn1k}) and (\ref{enu:fnorth}).

Let the sequence $\{q(n)\}$ above be fixed. The {}``almost orthogonality''
condition (\ref{enu:fnorth}) implies the {}``almost Pythagorean''
equality:\[
||f_{n-1}||^{2}+||f_{n}-f_{n-1}||^{2}=||f_{n}||^{2}+o(1).\]
From section \ref{sub:lnKnun} we have $||f_{n}||^{2}=g(l(n))$, so
$||f_{n}-f_{n-1}||^{2}=o(1)$. Together with (\ref{enu:defNn}) and
(\ref{enu:fnfn1k}) this easily gives\begin{equation}
||\widehat{f_{n}-f_{m}}||_{\infty }\rightarrow 0\textrm{ as }n,m\rightarrow \infty \label{eq:fnfm20}\end{equation}

\subsection{\label{sub:defmu}}

Let $\mu $ be the weak limit of the measures $\mu _{n}(dt):=u_{n}(t)dt$.
Clearly it is a positive measure supported on $K:=\cap K_{n}$ and
$|K|=0$. Define $F$ by (\ref{eq:poisson}) and $F_{n}$ by the same
formula with $\mu $ replaced by $\mu _{n}$. Then $F_{n}\rightarrow F$
uniformly on compacts inside the unit disc $\mathbb{D}$. Therefore
each coefficient $c(k)$ of the expansion (\ref{eq:defF}) may be
obtained as the limit of the corresponding coefficients $c_{n}(k)$
which are just $\widehat{f_{n}}(k)$, so (\ref{eq:fnfm20}) implies
(\ref{eq:cno1}) which finishes the proof.

\begin{rem*}
It should be noted that our use of Carleson's convergence theorem
to prove the equivalence of theorems \ref{thm:Mensh} and \ref{thm:complx}
is unnecessary, since theorem \ref{thm:Mensh} may be proved directly
using the fact (which is easy to see) that the function $f$ defined
in section \ref{sub:defFf} is smooth on any closed arc $J\subset K^{c}$.
\end{rem*}

\section{REMARKS}

\subsection{}

In contrast to Menshov's original example, the series (\ref{eq:defser})
in theorem \ref{thm:Mensh} can not be the Fourier series of a measure.
Indeed , if \[
\mu \sim \sum c(n)e^{int}\quad \sum _{n<0}|c(n)|^{2}<\infty ,\]
then $\mu $ must be absolutely continuous, and cannot generate a
null-series. See, for example, \cite[sect. VIII.12]{B64}.

\subsection{}

Another contrast with the {}``non-analytic'' situation appears when
one considers the size of the exceptional set. It is well known that
a null series (\ref{eq:defser}) may converge to zero outside a {}``thin''
compact (of zero Hausdorff dimension). On the other hand the following
proposition is true\begin{ialt}If a series (\ref{eq:sumanal}) converges
to $f\in L^{1}(\mathbb{T})$ everywhere on $\mathbb{T}$ outside some
set of dimension $<1$ then it is the Fourier series of $f$.\end{ialt}

This follows from a Phragm\'en-Lindel\"of type theorem for analytic
functions in $\mathbb{D}$ of slow growth, see \cite[theorem 5]{B92}.

\subsection{}

Let $\mathcal{P}$ be the class of functions in $L^{2}$ which can be represented
by an a.e.\ converging sum (\ref{eq:sumanal}).
Theorem \ref{thm:complx}
shows us that $\mathcal{P}\setminus H^{2}$ is non-trivial. Further,
the proof actually gives a little more: $(\mathcal{P}\setminus H^{2})\cap L^{\infty }\neq \emptyset $.
The class $\mathcal{P}$ has some interesting properties. We plan
to analyze it in a separate paper.

\end{document}